\documentclass[bj,authoryear]{imsart}

\usepackage{etex}
\reserveinserts{18}

\usepackage{amsmath,amsthm,amssymb,amsfonts}
\usepackage{mathtools}
\usepackage{parallel,accents}

\usepackage{amsxtra,amstext,latexsym,dsfont} 

\normalfont
\usepackage{bm}

\usepackage[comma]{natbib}

\usepackage[pdftex]{graphicx}
\usepackage[dvipsnames]{xcolor}

\usepackage{pgf,tikz}
\usetikzlibrary{snakes}
\usetikzlibrary{backgrounds}
\usetikzlibrary{arrows}
\usetikzlibrary{patterns}
\usetikzlibrary{fit}
\usetikzlibrary{calc}
\usetikzlibrary{matrix}

\usepackage{lscape} 

\usepackage{array,longtable,booktabs,float,rotating}
\usepackage{dcolumn,ctable}

\usepackage[labelfont={bf},ruled,labelsep=period,small]{caption}

\usepackage[flushleft]{threeparttable}

\usepackage[linkcolor=black,colorlinks=true,citecolor=black,
bookmarks=true,bookmarksopen=true,bookmarksnumbered=true,bookmarksopenlevel=0,%
hyperindex=true,pdfpagemode=UseOutlines,linktocpage=true,pdfpagelayout=TwoColumnLeft,
pdftitle=Biostatistics,pdfstartview=FitH,pdffitwindow=true,urlcolor=blue]{hyperref}

\usepackage{ragged2e}

\usepackage{multirow,multicol}

\usepackage{fancyhdr,extramarks}

 \makeatletter
 \newcommand\sub{\@startsection%
     {subsubsection}{5}{0mm}{-1\baselineskip}{.01\baselineskip}%
     {\normalfont\itshape}}
 \renewcommand\subsubsection{\@startsection%
     {subsubsection}{3}{0mm}{-1\baselineskip}{.01\baselineskip}%
     {\normalfont\itshape}}
 \makeatother

\makeatletter
        \newcommand\Appendix[2][?]{%
            \refstepcounter{section}%
            \addcontentsline{toc}{appendix}%
                {\protect\numberline{\appendixname~\thesection}#1}%
            {\raggedleft\bfseries \appendixname\
                \thesection\par \centering#2\par}%
                \sectionmark{#1}%
                \@afterheading
                \addvspace{\baselineskip}}
        \newcommand\sAppendix[1]{%
            \raggedleft\bfseries\appendixname\par
            \@afterheading\addvspace{\baselineskip}}
    \makeatother

\normalsize \makeindex

\newcolumntype{A}{>{\centering}p{100pt}}
\newlength\savedwidth

\def\coldot{.}%
{\catcode`\.=\active%
    \gdef.{$\egroup\setbox2=\hbox to \dimen0 \bgroup$\coldot}}
\def\rightdots#1{%
    \setbox0=\hbox{$1$}\dimen0=#1\wd0%
    \setbox0=\hbox{$\coldot$}\advance\dimen0 \wd0%
    \setbox2=\hbox to \dimen0 {}%
    \setbox0=\hbox\bgroup\mathcode`\.="8000 $}
\def\endrightdots{$\hfil\egroup\box0\box2}

\newcolumntype{d}[1]{D{.}{.}{#1}}
\newcolumntype{A}{>{\centering}p{100pt}}
\newcolumntype{.}{D{.}{.}{-1}}
\newcolumntype{P}[2]{>{#1\raggedright\arraybackslash}p{#2}}

\DeclareFontFamily{U}{euc}{}
\DeclareFontShape{U}{euc}{m}{n}{<-6>eurm5<6-8>eurm7<8->eurm10}{}%

\newtheorem{thm}{Theorem}
\theoremstyle{plain}      \newtheorem{lem}{Lemma}
\theoremstyle{plain}      \newtheorem{dfn}{Definition}
\theoremstyle{plain}      
\theoremstyle{plain}      
\theoremstyle{definition} \newtheorem{rmk}{Remark}
\theoremstyle{definition} 
\theoremstyle{definition} \newtheorem{exa}{Example}
\theoremstyle{plain} 
\theoremstyle{definition} 
\theoremstyle{plain} 
\theoremstyle{definition} 
\theoremstyle{definition} 
\theoremstyle{definition} 

\newcounter{nctr}
\usepackage[rightbars,color]{changebar}
\cbcolor{blue}

\usepackage{url}
\usepackage{fancyvrb}

\newcommand\ti{\textit}

\newcommand\bb{\mathbb}

\newcommand\te{\text}
\newcommand\ma[1]{\te{\bf{#1}}}
\newcommand\ca{\mathcal}

\newcommand\op{\operatorname}

\newcommand\as{^\ast}

\newcommand\arginf{\operatornamewithlimits{arginf}}

\newcommand\argmin{\operatornamewithlimits{argmin}}

\newcommand\E{\bb{E}}

\newcommand\for{\,\,\forall\,\,}

\newcommand\lb{\lbrace}
\newcommand\lt{\left}

\newcommand\p{\bb{P}}           

\newcommand\pri{^\prime}

\newcommand\rb{\rbrace}
\newcommand\rt{\right}

\newcommand\stack{\stackrel} 

\newcommand\tth{^\text{th}}

\newcommand\N{\bb{N}}  
\newcommand\Na{\bb{N}} 
\newcommand\R{\bb{R}}  

\newcommand\bX{\ma{X}}

\newcommand\cB{\ca{B}} 
\newcommand\cD{\ca{D}} 
\newcommand\cE{\ca{E}} 
\newcommand\cF{\ca{F}} %
\newcommand\cG{\ca{G}} %
\newcommand\cI{\ca{I}} 
\newcommand\cM{\ca{M}} 
\newcommand\cX{\ca{X}} %
\newcommand\sig{\sigma}

\newcommand\Th{\Theta}

\DeclareMathOperator*\climsup{limsup}
\DeclareMathOperator*\climinf{liminf}
\DeclareMathOperator*\limcsup{lim\,\oline{\op{sup}}}
\DeclareMathOperator*\Limsup{Limsup}

\newcommand*\oline[1]{%
  \vbox{%
    \hrule height 0.5pt
    \kern0.5ex
    \hbox{%
      \kern-0.1em
      \ifmmode#1\else\ensuremath{#1}\fi
      \kern-0.1em
    }
  }
}

\startlocaldefs
\endlocaldefs

\begin{document}
\sloppy

\begin{frontmatter}

\title{Strong Consistency of Fr\'{e}chet Sample
       Mean Sets for Graph-Valued Random Variables}
\runtitle{Consistency of Fr\'{e}chet Sample Mean Sets}

\author{\fnms{Cedric E.} \snm{Ginestet}\corref{}\ead[label=e1]{cedric.ginestet@kcl.ac.uk}\thanksref{t1}}
\thankstext{t1}{This work was supported by a grant from the Air Force
Office for Scientific Research (AFOSR), whose grant number is
FA9550-12-1-0102; and by a fellowship from the UK National Institute
for Health Research (NIHR) Biomedical Research Center for Mental
Health (BRC-MH) at the South London and Maudsley NHS Foundation Trust
and King's College London. Some portions of this work was
conducted during a visit of the author at the University of Warwick,
which was generously supported by the Center for Research in
Statistical Methodology (CRiSM). We would also like to thank Pietro
Rigo, Eric D. Kolaczyk, Thomas Nichols and Wilfrid S. Kendall for
useful discussion.} 
\address{Department of Mathematics and Statistics \\
Boston University, Boston, MA.}
\affiliation{Boston University}

\runauthor{Cedric E. Ginestet}
\begin{abstract}
    The Fr\'{e}chet mean or barycenter generalizes the idea of
    averaging in spaces where
    pairwise addition is not well-defined. In general metric spaces,
    the Fr\'{e}chet sample mean is
    not a consistent estimator of the theoretical Fr\'{e}chet mean. 
    For graph-valued random variables, for instance, the Fr\'{e}chet
    sample mean may fail to converge to a unique value. Hence, it
    becomes necessary to consider the convergence of sequences of sets
    of graphs. We show that a specific type of almost
    sure (a.s.) convergence for the Fr\'{e}chet sample mean previously
    introduced by \citet{Ziezold1977} is, in fact, equivalent to
    the Kuratowski outer limit of a sequence of Fr\'{e}chet sample means.
    Equipped with this outer limit, we provide a new proof of the
    strong consistency of the Fr\'{e}chet sample mean for graph-valued
    random variables in separable (pseudo-)metric space. Our proof
    strategy exploits the fact that the metric of interest is bounded,
    since we are considering graphs over a finite number of vertices.
    In this setting, we describe two strong laws of large numbers for both
    the restricted and unrestricted Fr\'{e}chet sample means of all
    orders, thereby generalizing a previous result, due to 
    \citet{Sverdrup1981}.
\end{abstract}
\end{frontmatter}

\section{Introduction}\label{sec:intro}
All statistics are summaries. The epitome of these
summaries is the sample mean, and its theoretical analog, the expected
value. In an inspired monograph, \citet{Frechet1948} generalized this
concept to any abstract metric space. He showed that the sole
requirement for the definition of a mean element is the specification
of a metric on the space of interest. Once this metric has been
chosen and a probability measure has been defined on that metric
space, the Fr\'{e}chet mean is simply the element that minimizes
the sum of the squared distances from all the
elements in that space. The Fr\'{e}chet mean generalizes other notions of
means in abstract spaces, such as the centroid in
Euclidean geometry, the barycenter or center of mass in physics, the
Procrustean mean in shape spaces \citep{Le1998}, and the Karcher mean
on Riemannian manifolds
\citep{Karcher1977}. The sample version of the Fr\'{e}chet mean 
can naturally be expressed using cumulative addition instead of the
expectation, thereby producing a convex combination operator on 
metric spaces with both negative and positive Alexandrov curvature
\citep{Ginestet2012b}. 

The object of this paper is to characterize the asymptotic behavior of
the Fr\'{e}chet sample mean in separable metric
spaces with a bounded metric. We are here especially interested in metric spaces of
simple graphs. Separability is a relatively mild topological assumption
likely to be satisfied in most applications. The boundedness of the
metric, however, is a more stringent condition. Nonetheless, there
is a range of modern statistical applications for which the metric of
interest is likely to be bounded. In bioinformatics, the use of the
\citet{Hamming1950} distance on finite alphabets, such as stretches of DNA for
instance, naturally gives rise to such assumptions
\citep{He2004}. Similarly, the comparison of families of
networks with a given number of nodes, as commonly done in neuroscience
\citep{Ginestet2011b} may similarly generate bounded
metric spaces; albeit the combinatorial nature of these metrics may
lead to bounds that increase factorially with the number of nodes in
these networks. 

The asymptotic properties of the Fr\'{e}chet sample mean have been studied
by several authors. \citet{Ziezold1977} proved a strong law of large
numbers for Fr\'{e}chet sample means defined in separable
pseudo-metric spaces, where the metric is not assumed to satisfy the coincidence axiom. 
This a.s.~ convergence result has also been demonstrated for compact metric
spaces by \citet{Sverdrup1981}. The perspectives adopted by these two
authors are very different in nature. Given the fact that \citet{Sverdrup1981}
does not cite the work of \citet{Ziezold1977}, and because the work of
the latter was published in a conference proceedings, it is probable
that \citet{Sverdrup1981} was not cognisant of Ziezold's proof
technique. 


The properties of sample Fr\'{e}chet means on Riemannian manifolds
have been particularly well-studied
\citep{Bhattacharya2002,Bhattacharya2005,Bhattacharya2012}. 
When the Fr\'{e}chet mean is assumed to be unique, the theorem of
\citet{Sverdrup1981} has been generalized by
\citet{Bhattacharya2003} for \textit{proper} metric spaces. Recall
that a metric space is proper, if and only if every bounded closed
subsets of that space is compact \citep{Sahib1998,Yang2011}. 
By the Hopf-Rinow theorem, every complete and connected Riemannian
manifold is a proper metric space. Thus, \citet{Bhattacharya2003} have
weakened the compactness assumption made by \citet{Sverdrup1981}, and
their strong law of large numbers apply to manifolds, under some very
mild conditions. Recently, \citet{Kendall2011} have further generalized
these results with a weak law of large numbers and a central limit
theorem for sequences of Fr\'{e}chet sample means based on non-iid random
variables taking values on a Riemannian manifold. 

Here, we consider sequences of random variables taking
values in separable pseudo-metric spaces with a bounded metric. 
Using boundedness, we provide a different proof of the strong
consistency of the Fr\'echet sample mean from the one of 
\citet{Ziezold1977}. In addition, we generalize the results of 
 \citet{Sverdrup1981} on \textit{restricted} Fr\'echet sample means.
The restricted Fr\'{e}chet sample mean is the most `typical' quantity
chosen from the available sampled values. The computation of the
unrestricted Fr\'{e}chet sample mean in arbitrary metric spaces can 
indeed prove to be arduous, since this necessarily requires a minimization
over a complex space. The difficulties that arise when estimating the
Fr\'{e}chet mean in shape spaces, for instance, have received special
attention \citep{Dryden1998,Kume2000,Le2001,Le2004}. Estimation issues
have also been addressed in spaces of covariance matrices, where a range of
different metrics can be considered
\citep{Arsigny2007,Dryden2009,Yang2011a}. For graph-valued random
variables, several metrics have been proposed in the literature, which
are NP-hard to minimize. The restricted Fr\'{e}chet mean may therefore
be useful in practice, as it greatly simplifies the minimization
procedure, by simply selecting the most typical element in the sample.

Importantly, we also clarify previous results on the asymptotic
consistency of the Fr\'{e}chet sample mean, by showing
that the modes of convergence studied by \citet{Ziezold1977} and
\citet{Sverdrup1981} are, in fact, equivalent to the consideration of
the Kuratowski outer limit of a
sequence of Fr\'{e}chet sample means. 
One of the core difficulties with the consideration of the asymptotic
properties of Fr\'{e}chet sample means is that such functions can be
\textit{multivalued}. That is, when the Fr\'{e}chet sample mean is not unique, we
obtain a random variable that is a set-valued function, which takes
values in the power set of $\cX$, or more precisely in the Borel
$\sigma$-algebra of $\cX$. It then
becomes necessary to consider the convergence of multivalued
functions. To this end, we resort to the tools of set-valued analysis,
as described by \citet{Aubin1990}. This difficulty leads us to
consider different `types' of convergence, depending on whether we require the
Fr\'{e}chet sample mean to converge, or are simply interested in
evaluating the asymptotic behavior of the outer limit of that
sequence \citep[see][for an introduction to set-valued random
variables]{Molchanov2005}.

The main innovation in this paper is our formal set-valued
perspective. Note that our approach differs from the one of 
\citet{Bhattacharya2012}, since we have allowed the metric spaces of
interest to be non-compact, and not necessarily equipped with a
manifold structure. In particular, we identify the key role
played by the Kuratowski outer limit when studying sequences of 
Fr\'{e}chet sample means. This paper therefore constitutes an
extension of the work of \citet{Ziezold1977} and \citet{Sverdrup1981}
to Fr\'echet means of all orders, and to restricted Fr\'echet means. 
Moreover, we have emphasized the importance of point functions and of the
Glivenko-Cantelli lemma. 

This paper is organized as follows. Firstly, we motivate this work with
a counterintuitive example of a graph-valued mean set that includes its
sample as a proper subset. This justifies our emphasis on set-valued
convergence throughout the rest of the paper. 
In section \ref{sec:sequence}, we then
introduce and study different types of a.s.~ convergence for sequences
of Fr\'{e}chet sample mean sets, and show through counterexamples why the
Kuratowski outer limit is adequate for this purpose.
In section \ref{sec:main result}, we prove the strong consistency
of the Fr\'{e}chet sample mean sets in bounded metric spaces. 
Finally, section \ref{sec:restricted} is devoted to
the description of the restricted versions of the Fr\'{e}chet sample
mean, and a generalization of a result due to \citet{Sverdrup1981} to
bounded metric spaces, for random variables with closed support. 
\begin{figure}[t]
    \centering
         \begin{tikzpicture}[font=\small,scale=1]
          \draw(0,-.5)node{$G_{1}$};
          \draw(.75,0)node(v0){}; \draw(-.75,0)node(v1){};
          \draw(1.25,1.5)node(v2){}; \draw(-1.25,1.5)node(v3){};
          \draw(0,2.25)node(v4){};
          \fill[fill=black](v0) circle (2.0pt); \fill[fill=black](v1)
          circle (2.0pt); \fill[fill=black](v2) circle (2.0pt);
          \fill[fill=black](v3) circle (2.0pt); \fill[fill=black](v4)
          circle (2.0pt);
          \foreach \x in {v1,v2} \draw (v0) -- (\x); \foreach \x in
          {v2,v3,v4} \draw (v1) -- (\x); \foreach \x in {v4} \draw
          (v2) -- (\x); \draw(v3) -- (v4);
        \end{tikzpicture}
        \hspace{1cm}
        \begin{tikzpicture}[font=\small,scale=1]
          \draw(0,-.5)node{$G_{2}$};
          \draw(.75,0)node(v0){}; \draw(-.75,0)node(v1){};
          \draw(1.25,1.5)node(v2){}; \draw(-1.25,1.5)node(v3){};
          \draw(0,2.25)node(v4){};
          \fill[fill=black](v0) circle (2.0pt); \fill[fill=black](v1)
          circle (2.0pt); \fill[fill=black](v2) circle (2.0pt);
          \fill[fill=black](v3) circle (2.0pt); \fill[fill=black](v4)
          circle (2.0pt);
          \foreach \x in {v1,v3,v4} \draw (v0) -- (\x); \foreach \x in
          {v3,v4} \draw (v1) -- (\x); \foreach \x in {v3} \draw (v2)
          -- (\x); \foreach \x in {v4} \draw (v3) -- (\x);
        \end{tikzpicture}        
        \hspace{1cm}
        \begin{tikzpicture}[font=\small,scale=1]
          \draw(0,-.5)node{$G_{3}$};
          \draw(.75,0)node(v0){}; \draw(-.75,0)node(v1){};
          \draw(1.25,1.5)node(v2){}; \draw(-1.25,1.5)node(v3){};
          \draw(0,2.25)node(v4){};
          \fill[fill=black](v0) circle (2.0pt); \fill[fill=black](v1)
          circle (2.0pt); \fill[fill=black](v2) circle (2.0pt);
          \fill[fill=black](v3) circle (2.0pt); \fill[fill=black](v4)
          circle (2.0pt);
          \foreach \x in {v1,v3,v4} \draw (v2) -- (\x); \foreach \x in
          {v3,v4} \draw (v1) -- (\x); \foreach \x in {v3} \draw (v2)
          -- (\x); \foreach \x in {v4} \draw (v3) -- (\x);
        \end{tikzpicture}        
   \caption{A sample of graphs, $G_{i}=(V,E)$, over five vertices, denoted by
            $G_{i}\in\cG_{5}$. 
            \label{fig:motivation}}
\end{figure}
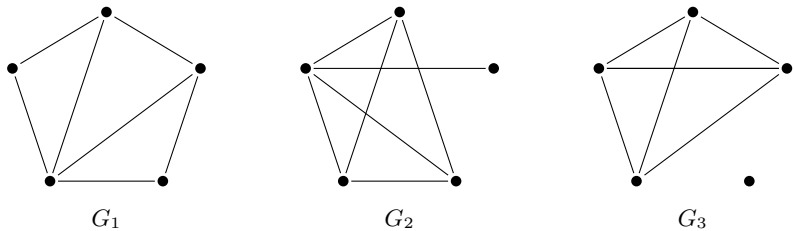

\section{Motivating Example: Graph Means}\label{sec:motivation}
We are here especially interested in spaces of simple graphs,
$G_{i}:=(V,E)$ with $i=1,\ldots,n$, which have a fixed
number of vertices, $N_{v}:=|V(G_{i})|$, but their edge set, $E(G_{i})$
may vary. A graph is said to be simple, when it does not contain multiple
edges, loops or weighted edges. Throughout this paper, we will assume that there
exists a probability measure on the space of all such simple graphs.
A sample of three such simple graphs for $N_{v}=7$ is
given in figure \ref{fig:motivation}. 

Statistically, one may be interested in computing the mean graph for
this type of random variables. Such a mean quantity can be defined as the Fr\'echet
mean of that variable with respect to some distance function on the
space of interest. A standard distance function on spaces of graphs is the
Hamming distance, which is defined as follows for any two graphs $G=(V,E)$
and $G\pri=(V,E\pri)$ with $N_{v}$ vertices, 
\begin{equation}
   d_{H}(G,G\pri) := \sum_{i<j} \cI\{e_{ij}\neq e\pri_{ij}\}.
   \notag
\end{equation}
We denote by $\cG_{N_{v}}$ the space of all simple graphs with $N_{v}$ vertices.
Given a graph-valued random variable on $\cG_{N_{v}}$, the mean value for a sample of
$n$ realizations is then given by the element in $\cG_{N_{v}}$, which minimizes the
squared distances to all the graphs in the sample considered. For
general graph-valued random variables, however, such a mean element needs not
be unique. 

In figure \ref{fig:hamming}, we consider a sample of $n=2$ graphs
$S_{1}$ and $S_{2}$ with $N_{v}=4$ vertices. Using the Hamming distance, the
Fr\'echet mean graphs are the following elements of $\cG_{4}$,
\begin{equation}
    \Th := \argmin_{G\pri\in\cG_{4}} \sum_{l=1}^{n}\sum_{i<j}
    \cI\{e^{(l)}_{ij}\neq e\pri_{ij}\}.
    \notag
\end{equation}
One can easily verify that the Fr\'echet mean is given
by a set of four different simple graphs, as shown in figure \ref{fig:hamming}.
Hence, in this setting, we obtain the paradoxical result that the
sample is a \textit{proper} subset of the mean. This is somewhat
counterintuitive, since we generally expect an average value to
summarize information, and therefore to be more `concentrated' than
the sample values on which the mean is based. 

Observe that the Hamming distance is here a bounded metric. 
In the sequel, we will consider the more general case of random
variables taking values in separable metric spaces with bounded metrics,
which encompasses graph-valued random variables, as a special case.
Other popular choices of distance functions include the graph edit
distance \citep{Gao2010}, and maximum common subgraph distance
\citep{Bunke1997}.
\begin{figure}[t]
      \centering
      \textit{(a) Sample of simple graphs, $S$, with $n=2$.}\\
      \vspace{.3cm}
      \begin{tikzpicture}[font=\small,scale=1]
        \draw(.5,-.5)node{$S_{1}$};
        \draw(0,0)node(v1){}; \draw(1,0)node(v2){};
        \draw(1,1)node(v3){}; \draw(0,1)node(v4){};
        \fill[fill=black](v1) circle (2.0pt); \fill[fill=black](v2)
        circle (2.0pt); \fill[fill=black](v3) circle (2.0pt);
        \fill[fill=black](v4) circle (2.0pt);
        \draw(v1)--(v2); \draw(v2)--(v3); \draw(v3)--(v4);
      \end{tikzpicture}
      \hspace{1cm}
      \begin{tikzpicture}[font=\small,scale=1]
        \draw(.5,-.5)node{$S_{2}$};
        \draw(0,0)node(v1){}; \draw(1,0)node(v2){};
        \draw(1,1)node(v3){}; \draw(0,1)node(v4){};
        \fill[fill=black](v1) circle (2.0pt); \fill[fill=black](v2)
        circle (2.0pt); \fill[fill=black](v3) circle (2.0pt);
        \fill[fill=black](v4) circle (2.0pt);
        \draw(v1)--(v2); \draw(v3)--(v4); \draw(v4)--(v1);
      \end{tikzpicture}      
      \\
      \vspace{.5cm}     
      \textit{(b) Fr\'echet Mean, $\Th$, for this sample.}\\
      \vspace{.3cm}
      \begin{tikzpicture}[font=\small,scale=1]
        \draw(.5,-.5)node{$\Theta_{1}$};
        \draw(0,0)node(v1){}; \draw(1,0)node(v2){};
        \draw(1,1)node(v3){}; \draw(0,1)node(v4){};
        \fill[fill=black](v1) circle (2.0pt); \fill[fill=black](v2)
        circle (2.0pt); \fill[fill=black](v3) circle (2.0pt);
        \fill[fill=black](v4) circle (2.0pt);
        \draw(v1)--(v2); \draw(v2)--(v3); \draw(v3)--(v4);
      \end{tikzpicture}
      \hspace{1cm}
      \begin{tikzpicture}[font=\small,scale=1]
        \draw(.5,-.5)node{$\Theta_{2}$};
        \draw(0,0)node(v1){}; \draw(1,0)node(v2){};
        \draw(1,1)node(v3){}; \draw(0,1)node(v4){};
        \fill[fill=black](v1) circle (2.0pt); \fill[fill=black](v2)
        circle (2.0pt); \fill[fill=black](v3) circle (2.0pt);
        \fill[fill=black](v4) circle (2.0pt);
        \draw(v1)--(v2); \draw(v2)--(v3); \draw(v3)--(v4); \draw(v1)--(v4);
      \end{tikzpicture}
      \hspace{1cm}
      \begin{tikzpicture}[font=\small,scale=1]
        \draw(.5,-.5)node{$\Theta_{3}$};
        \draw(0,0)node(v1){}; \draw(1,0)node(v2){};
        \draw(1,1)node(v3){}; \draw(0,1)node(v4){};
        \fill[fill=black](v1) circle (2.0pt); \fill[fill=black](v2)
        circle (2.0pt); \fill[fill=black](v3) circle (2.0pt);
        \fill[fill=black](v4) circle (2.0pt);
        \draw(v1)--(v2); \draw(v3)--(v4);
      \end{tikzpicture}
      \hspace{1cm}
      \begin{tikzpicture}[font=\small,scale=1]
        \draw(.5,-.5)node{$\Theta_{4}$};
        \draw(0,0)node(v1){}; \draw(1,0)node(v2){};
        \draw(1,1)node(v3){}; \draw(0,1)node(v4){};
        \fill[fill=black](v1) circle (2.0pt); \fill[fill=black](v2)
        circle (2.0pt); \fill[fill=black](v3) circle (2.0pt);
        \fill[fill=black](v4) circle (2.0pt);
        \draw(v1)--(v2); \draw(v3)--(v4); \draw(v4)--(v1);
      \end{tikzpicture}      
    \caption{The sample of graphs in (a) is here a proper subset of the
        graph mean in (b), such that $S\subset\Th$, where the Fr\'echet mean,
        $\Th$ is computed with respect to the Hamming distance on the
        space of all simple graphs with $N_{v}=4$ vertices.
        \label{fig:hamming}}
\end{figure}

\section{Sequences of Fr\'{e}chet Sample Means}\label{sec:sequence}
\subsection{Empirical and Theoretical Fr\'{e}chet Means}
A separable space $\cX$ is endowed with a metric $d:\cX\times\cX\mapsto\bb{R}^{+}$. 
This produces a metric space, $(\cX,d)$, with elements $x$. 
Let a probability space be denoted by $(\Omega,\cF,\p)$, and define a random
variable, $X$, on that space, which takes values in $(\cX,\cB)$.
Here, $\cB$ is the Borel $\sigma$-algebra generated by the topology,
$\tau$ on $\cX$, induced by $d$. The triple $(\Omega,\cF,\p)$ is
assumed to be \textit{complete}, in the sense that every subset of
every null set is measurable. This
is particularly convenient for constructing product spaces based on
$\Omega$ that remain well-behaved. In addition, we define 
$\mu(B):=(\p\circ X^{-1})(B)$, for every $B\in\cB$. Naturally, $X$ is
here assumed to be $(\cF,\cB)$-measurable. 
Such a random variable will be termed an \textit{abstract-valued} random
variable, which will be contrasted with the more standard
real-valued random variables. 

In this setting, we compute the most `\ti{central}' element. This is
the element that has the smallest expected distance to all other
elements in $\cX$. This approach allows us to define the following
moments \citep{Frechet1948},
\begin{equation}
     \Theta^{r} := \arginf_{x\pri\in\cX} \int\limits_{\cX}d(x,x\pri)^{r}
     d\mu(x), 
     \quad\te{and}\quad
     \sig^{r} := \inf_{x\pri\in\cX}\int\limits_{\cX}d(x,x\pri)^{r}d\mu(x),
     \label{eq:frechet theory}
\end{equation}
for every $0<r<\infty$, and where $\Theta^{r}\subseteq\cX$. Observe that
we are using the superscript $r$ on the Fr\'{e}chet variance as a simple
marker of the order of the exponentiated metric. Thus, in general, it
will not be true that $(\sig^{r})^{1/r}$ simplifies to $\sig^{1}$.

These are commonly referred to as the Fr\'{e}chet
mean and variance when $r=2$. For other choices of $r$, we will refer
to these different Fr\'{e}chet moments as Fr\'{e}chet moments of
\textit{order} $r$. Note that if
the \textit{infimum} of $\E[d(x,x\pri)^{r}]$ exists, then it is
unique. However, the \textit{argument of the infimum} may not
necessarily exist and may not be unique. If such an argument does not
exist, then $\Theta^{r}=\varnothing$. When the minimizer is not
unique, the ensemble of minimizers is
sometimes referred to as the \ti{Fr\'{e}chet mean set}. In particular, observe that
if $\Theta$ is not a singleton, $\sig^{2}=\E[d(X,\theta)^{2}]$
for any $\theta\in\Theta$, will not, in general, be equivalent to
$\E[d(X,\Theta)^{2}]$, where the distance between an element $x$
and a non-empty subset $A$ of $\cX$ is defined as
$d(x,A):=\inf\{d(x,y):y\in A\}$, with $d(x,\varnothing)=\infty$.
In this paper, Fr\'{e}chet mean and Fr\'{e}chet mean set will be used
interchangeably. Observe that when $\cX$ is a Hilbert space, endowed
with the inner product metric, then there exists a unique global
minimizer and $\Theta$ is therefore a singleton. 

Analogously, for a given sequence of abstract-valued random variables
$X_{i}:\Omega\mapsto\cX$, for every $i=1,\ldots,n$, one may define the
following Fr\'{e}chet sample moments of the $r\tth$ order 
\begin{equation}
      \widehat{\Theta}^{r}_{n}:=\arginf_{x\pri\in\cX}
      \frac{1}{n}\sum_{i=1}^{n}d(X_{i},x\pri)^{r} 
     \quad\te{and}\quad
     \widehat{\sig}_{n}^{r} := \inf_{x\pri\in\cX}\frac{1}{n}\sum_{i=1}^{n}d(X_{i},x\pri)^{r}.
     \label{eq:frechet sample}
\end{equation}
Observe that, even for the sample versions of the Fr\'{e}chet moments,
these infima meed not be attained, and therefore these quantities may
be empty for each $n$. When there is no ambiguity as to the order of
$\widehat{\Theta}^{r}_{n}$,
we will simply refer to this quantity as $\widehat{\Theta}_{n}$, and
similarly for $\Theta$. In the sequel, an element of $\Theta$ and
an element of $\widehat{\Theta}_{n}$ will be respectively denoted
by $\theta$ and $\hat{\theta}_{n}$. Our interest will mainly lie in
considering Fr\'{e}chet moments of the second order, albeit some
examples will also be studied where $r=1$. It is easy to see that the
Fr\'{e}chet mean and Fr\'{e}chet sample mean are closed subsets of 
$\cX$, if $\cX$ is Polish. 
\begin{lem}\label{lem:closed}
    For any space $(\cX,d)$, $\Theta^{r}$ and
    the $\,\widehat{\Theta}^{r}_{n}$'s are closed in $\cX$, for every $r\geq1$.
\end{lem}
\begin{proof}
    Clearly, if $\Theta^{r}=\varnothing$, then
    $\op{cl}(\Theta^{r})=\Theta^{r}$ and similarly for the $\widehat{\Theta}^{r}_{n}$'s.
    Now, fix $r=1$, and consider the Fr\'{e}chet mean set
    $\Theta\subseteq\cX$. Recall that the boundary of $\Theta$ is defined as 
    $\partial(\Theta):=\big\{ x\in\cX: d(\Theta,x)=d(\Theta^{C},x)=0\big\}$,
    where $\Theta^{C}:=\cX\setminus\Theta$. We proceed by
    contradiction. Assume that $\theta_{0}\in\partial(\Theta)$ and
    $\theta_{0}\notin\Theta$, then it follows that there exists
    $\theta\in\Theta$, such that by the triangle inequality,
    $d(\theta_{0},X)\leq d(\theta_{0},\theta) + d(\theta,X)$, for
    every $X\in\cX$. Taking the expectation, this gives
    \begin{equation}
         \E[d(\theta_{0},X)]\leq d(\theta_{0},\theta) +
         \E[d(\theta,X)] = \inf_{x\pri\in\cX}\E[d(X,x\pri)],
         \notag
    \end{equation}
    since $d(\theta_{0},\Theta)=0$, and using the definition of
    $\Theta$ in equation (\ref{eq:frechet theory}). Thus, $\theta_{0}$
    is optimal with respect to the infimum over $\cX$. However, we have
    assumed that $\theta_{0}\notin\Theta$, which leads to a contradiction,
    and therefore $\partial(\Theta)\subseteq\Theta$.

    Next, consider the case of $r>1$. Through a classical result on
    metric spaces \citep[see, for instance][p.229]{Frechet1948}, we have 
    \begin{equation}
         \Big(\E[d(\theta_{0},X)^{r}]\Big)^{1/r} 
         \leq \Big(\E[d(\theta_{0},\theta)^{r}]\Big)^{1/r} 
         + \Big(\E[d(\theta,X)^{r}]\Big)^{1/r},
         \notag
    \end{equation}
    for every $r>1$, and the result immediately follows, using the
    same argument. The proof is identical for the
    $\widehat{\Theta}^{r}_{n}$'s.
\end{proof}

\subsection{Convergence of Fr\'{e}chet Sample Mean Sets}\label{sec:as}
In this section, we study and compare different modes of convergence
for set-valued random variables. In particular, note that our chosen modes of
convergence differ from the ones used by \citet{Bhattacharya2012},
since we are not here assuming the compactness of the underlying
metric space $\cX$. Moreover, the target Fr\'echet mean set is also
allowed to be empty, thereby making it difficult to implement the
methods of \citet{Bhattacharya2012}. 

For the Fr\'{e}chet sample mean and its theoretical analogue, a.s.~
convergence could be defined in $(\cX,d)$ using sequences of random
sets as follows,
\begin{equation}
   \p\lt[ \lt\lb \omega\in\Omega: 
   \widehat{\Theta}_{n}(\omega) \to \Theta \rt\rb\rt] = 1,
   \label{eq:as definition1}
\end{equation}
where observe that $\Theta$ is here treated as a fixed subset of $\cX$.
The event in equation (\ref{eq:as definition1}) will have probability
one if the sequence of random sets, denoted $\widehat{\Theta}_{n}$, converges
a.s.~ in a set-theoretical sense such that
\begin{equation}
    \climinf_{n\to\infty}\widehat{\Theta}_{n}(\omega) =
    \climsup_{n\to\infty}\widehat{\Theta}_{n}(\omega) = \Theta,
    \label{eq:convergence criterion}
\end{equation}
for almost every $\omega\in\Omega$, and where $\op{liminf}
S_{n}:=\bigcup^{\infty}_{n=1}\bigcap^{\infty}_{m=n}S_{n}$, 
and $\op{limsup}S_{n}:=\bigcap^{\infty}_{n=1}\bigcup^{\infty}_{m=n}S_{n}$
denote the standard inner and outer
limits of a sequence of subsets of $\cX$. For most
purposes, however, this type of convergence is too strong. In fact, this criterion 
does not hold for Fr\'{e}chet sample means defined with respect to
general abstract-valued random variables. There are many non-trivial
examples of sequences of 
Fr\'{e}chet sample means that diverge. Consider the following example adapted from
the three-dimensional case described by \citet{Sverdrup1981}. 
\begin{figure}[t]
\centering
\tikzstyle{background rectangle}=[draw=gray!30,fill=gray!5,rounded corners=1ex]
\begin{tikzpicture}[font=\small,scale=.9,show background rectangle]
  \draw(0,3)node{\textbf{(a)}};  
  \draw(-2.25,0) -- (2,0);
  \draw(2,0)node[anchor=west]{$\cX$};
  \draw(-2,0) -- (-2,2.5);
  \draw(-2.0-.1,1) -- (-2.0+.1,1);
  \draw(-2.0-.1,2) -- (-2.0+.1,2);
  \draw(-2.0-.1,1)node[anchor=east]{$1$};
  \draw(-2.0-.1,2)node[anchor=east]{$2$};
  \draw(-1,-.1) -- (-1,.1);
  \draw(0,-.1) -- (0,.1);
  \draw(1,-.1) -- (1,.1);
  \draw(-1,-.1)node[anchor=north]{$-1$};
  \draw(0,-.1)node[anchor=north]{$0$};
  \draw(1,-.1)node[anchor=north]{$1$};
  \draw[-*](-1,0) -- (-1,1/2+.1);
  \draw[-*](1,0) -- (1,1/2+.1);
  \draw(-2,1/2-.05)node[anchor=east]{$\p[X=x]$};
  \draw[dashed](-2,1/2) -- (2,1/2);
  \draw(1,1)node[anchor=south]{$\E[d(X,x\pri)^{1}]$};
  \draw[very thick, color=blue!50](-1,1) -- (1,1);
\end{tikzpicture}
\begin{tikzpicture}[font=\small,scale=.9,show background rectangle]
  \draw(0,3)node{\textbf{(b)}};  
  \draw(-2.25,0) -- (2,0);
  \draw(2,0)node[anchor=west]{$\cX$};
  \draw(-2,0) -- (-2,2.5);
  \draw(-2.0-.1,1) -- (-2.0+.1,1);
  \draw(-2.0-.1,2) -- (-2.0+.1,2);
  \draw(-2.0-.1,1)node[anchor=east]{$1$};
  \draw(-2.0-.1,2)node[anchor=east]{$2$};
  \draw(-1,-.1) -- (-1,.1);
  \draw(0,-.1) -- (0,.1);
  \draw(1,-.1) -- (1,.1);
  \draw(-1,-.1)node[anchor=north]{$-1$};
  \draw(0,-.1)node[anchor=north]{$0$};
  \draw(1,-.1)node[anchor=north]{$1$};
  \draw[-*](-1,0) -- (-1,1/2+.1);
  \draw[-*](1,0) -- (1,1/2+.1);
  \draw(-2,1/2-.05)node[anchor=east]{$\p[X=x]$};
  \draw[dashed](-2,1/2) -- (2,1/2);
  \draw(1,2)node[anchor=south]{$\E[d(X,x\pri)^{2}]$};
  \draw[very thick, color=blue!50](0,1) -- (-1,2);
  \draw[very thick, color=blue!50](0,1) -- (1,2);
\end{tikzpicture}
  \caption{Metric and measure spaces considered in examples \ref{exa:divergent}
    and \ref{exa:closure}. In both panels, the closed interval $[-1,1]$
    is equipped with the Manhattan (or taxicab) metric, and two point
    masses are specified at $-1$ and $1$. Different
    Fr\'{e}chean inferences are conducted by taking $r=1$ and $r=2$ in panels (a)
    and (b), respectively. In the first case, the theoretical Fr\'{e}chet mean
    coincides with the median of $X$, whereas in panel (b), the
    theoretical Fr\'{e}chet mean coincides with the arithmetic
    mean. However, the sequence of Fr\'{e}chet \textit{sample} means
    diverge in both cases, when convergence is evaluated using
    set-valued liminf and limsup, as described in equation (\ref{eq:convergence
      criterion}).}
    \label{fig:examples}
\end{figure}
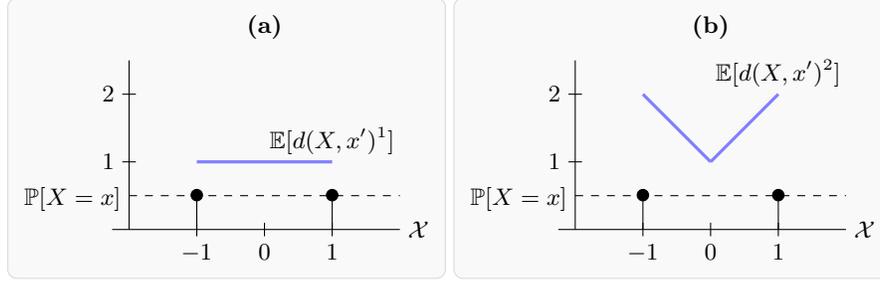
\begin{exa}\label{exa:divergent}
    Let the interval, $\cX:=[-1,1]\subset\R$, and equip this set with the
    usual Manhattan distance, defined as $d(x,y):=|x-y|$ for
    every $x,y\in\cX$. Additionally, let the random variable $X$, which
    takes values in $\cX$, and which satisfies the following
    $\p\lt[X=-1\rt]=\p\lt[X=1\rt]= 1/2$. This construction is illustrated
    in panel (a) of figure \ref{fig:examples}.
    The theoretical Fr\'{e}chet mean of order $r=1$ can be readily found as
    \begin{equation}
       \Theta^{1}=\arginf_{x\pri\in\cX} \sum_{x\in\{-1,1\}}d(x,x\pri)\p[x] = \cX,
       \notag
    \end{equation}
    since the energy function satisfies $\cE(x\pri):=\sum d(x,x\pri)\p[x]=1$ for every
    $x\pri\in\cX$. Here, the Fr\'{e}chet mean defined with respect to the Manhattan distance
    coincides with the \textit{median} of the real-valued random
    variable $X$ \citep{Feldman1966}.

    For the empirical Fr\'{e}chet mean, $\widehat{\Theta}^{1}_{n}$, first compute
    $S_{n}:=\sum_{i=1}^{n}X_{i}$. Clearly, the $S_{n}$'s are integer-valued.
    Observe the correspondence between the values of $S_{n}$ and the values taken
    by the Fr\'{e}chet sample mean. If the event $\{S_{n}=0\}$
    occurs, then it can easily be seen that
    $\widehat{\Theta}_{n}$ is equal to $\cX$. Similarly, 
    $\{S_{n}\geq 1\}$, and $\{S_{n}\leq -1\}$ respectively imply that 
    $\hat{\theta}_{n}=1$ and $\hat{\theta}_{n}=-1$. Now, 
    \begin{equation}
         \p\lt[\{S_{2n}=0\}\rt] = \binom{2n}{n}\lt(\frac{1}{2}\rt)^{2n}
         \approx (n\pi)^{-1/2},
         \notag
    \end{equation}
    for every $n$, using Stirling's approximation. Since
    $\p\lt[\{S_{n}=0\}\rt]$ is null, when $n$ is odd, it follows that 
    $\sum_{n=1}^{\infty} \p\lt[\{S_{n}=0\}\rt] < \infty$,
    and therefore by the Borel-Cantelli lemma, we have
    $\p\lt[\{S_{n}=0\}\;\te{i.o.}\rt]=0$, where i.o.~ means
    infinitely often. This implies that 
    $\p[\{\,\widehat{\Theta}_{n}=\cX\}\;\te{i.o.}]=0$, and hence 
    $\op{limsup}\widehat{\Theta}_{n}\neq\cX$.

    By using a similar argument, one can observe that 
    $\p[\{\,S_{n}\leq -1\}\;\te{i.o.}]=\p[\{\,S_{n}\geq 1\}\;\te{i.o.}]=1$, which implies that 
    $\p[\{\,\hat{\theta}_{n}=-1\}\;\te{i.o.}]=\p[\{\,\hat{\theta}_{n}=1\}\;\te{i.o.}]=1$, and
    therefore $\{-1,1\}$ is the limit superior of the sequence of
    Fr\'{e}chet mean sets. By contrast, there does not exist an $N>0$, such that 
    $\hat{\theta}_{n}=1$, for every $n\geq N$. An identical statement holds for 
    $\hat{\theta}_{n}=-1$, and therefore the limit inferior of
    $\widehat{\Theta}_{n}$ is empty. Thus,
    \begin{equation}
        \climsup_{n\to\infty}\widehat{\Theta}_{n}(\omega) = \{-1,1\} \supset 
        \climinf_{n\to\infty}\widehat{\Theta}_{n}(\omega) = \varnothing,
        \notag
    \end{equation}
    and the sequence of Fr\'{e}chet sample means diverges,
    as criterion (\ref{eq:convergence criterion}) is not
    satisfied. 
\end{exa}
The preceding example highlights two important aspects of the
asymptotic behavior of the Fr\'{e}chet sample mean set.
Firstly, the Fr\'{e}chet sample mean will in general fail to converge
in the sense that its outer and inner limits need not be
identical. In such cases, the sequence of Fr\'{e}chet sample means
exhibit an oscillatory property \citep[see][]{Feldman1966}.
Secondly, the limit superior of a sequence of Fr\'{e}chet
sample means may solely represent a subset of the theoretical Fr\'{e}chet
mean. Taken together, these two problems necessitate (i) the study of
the asymptotic behavior of the \textit{outer limit} of the
$\widehat{\Theta}_{n}$'s, and (ii) the consideration of the
convergence of the Fr\'{e}chet sample mean in terms of
\textit{set inclusion}, as a subset of the theoretical Fr\'{e}chet
mean. The passage from equations to inclusions is a natural step
in the generalization of singleton-valued analysis to 
set-valued analysis.

Example \ref{exa:divergent} leads to the formulation of a weaker type
of convergence, which can be expressed as the probability of the
following event, 
\begin{equation}
    \lt\lb \omega\in\Omega: \climsup_{n\to\infty}
    \widehat{\Theta}_{n}(\omega) \subseteq \Theta \rt\rb.
    \label{eq:as definition2}
\end{equation}
However, we here encounter a slightly different problem than the one
highlighted in our first example. This second issue can be
illustrated through another counterexample, which shows that this particular
type of a.s.~ convergence does not agree
with the analogous real-valued a.s.~ convergence. That is, the
reformulation of a given real-valued random variable into an
abstract-valued setting, equipped with the same topology produces a
divergent Fr\'{e}chet sample mean in terms of equation (\ref{eq:as
  definition2}). As a result, we obtain the somewhat
counterintuitive result that the arithmetic sample mean differs
from the corresponding Fr\'{e}chet sample mean. 
\begin{exa}\label{exa:closure}
    Consider the same setting described in example
    \ref{exa:divergent}, where now $r=2$ (see panel (b) of figure
    \ref{fig:examples}). One can immediately see that
    the theoretical Fr\'{e}chet mean is a singleton set, 
    \begin{equation}
       \Theta^{2}=\arginf_{x\pri\in\cX} \sum_{x\in\{-1,1\}}d(x,x\pri)^{2}\p[x] = 0,
       \notag
    \end{equation}
    which coincides with the expected value of the real-valued random
    variable $X$. For the Fr\'{e}chet sample mean, we know from example
    \ref{exa:divergent} that $\p[\{S_{n}=0\}\;\te{i.o.}]=0$ and
    therefore the probability of the sequence of empirical Fr\'{e}chet
    means including $\E[X]$ infinitely often is null. That is, for
    $r=2$, we have  $\p[\{\hat{\theta}_{n}=0\}\;\te{i.o.}]=0$. 
    Observe that the same is true for any other specific sequence
    of realizations of $X$. Consider the case of $S_{3n}=nx_{1} + 2nx_{2}$, 
    where $x_{1}=-1$ and $x_{2}=1$. For this subsequence, there exists
    a unique infimum, which is $\hat{\theta}_{n}=1/3$. The
    probability of this event occurring is as follows,
    \begin{equation}
        \p\lt[\{S_{3n}=nx_{1}+2nx_{2}\}\rt] = 
        \binom{3n}{n}\lt(\frac{1}{2}\rt)^{3n} \approx 
        (1/2)^{5n},
        \notag               
    \end{equation}
    which was approximated using Stirling's formulae. Clearly, all
    possible values of the Fr\'{e}chet sample mean of $X$ can be represented as a formulae of
    the form $nx_{1}+\alpha nx_{2}$, for some $\alpha\in\N$.
    Using the Borel-Cantelli lemma, it therefore follows that 
    there does not exist a point in $[-1,1]$ that $\hat{\theta}_{n}$
    will visit infinitely often, and hence
    $\op{limsup}\widehat{\Theta}_{n}=\op{liminf}\widehat{\Theta}_{n}=\varnothing$.
    By contrast, the arithmetic sample mean,
    $\bar{X}_{n}:=n^{-1}\sum_{i=1}^{n}X_{i}$ trivially converges to the
    expected value of $X$ a.s., since for every $\epsilon>0$, there exists
    an $N>1$, for which $d(\bar{X}_{n}(\omega),\E[X])<\epsilon$, for
    every $n\geq N$, for almost every $\omega\in\Omega$. Thus, for
    this example, we reach the counterintuitive conclusion that 
    $\bar{X}_{n}\notin \op{limsup}\widehat{\Theta}_{n}$, for every $n$.
\end{exa}
This paradoxical disagreement
between the divergence of the Fr\'{e}chet sample mean and the
classical convergence of the arithmetic sample mean in such a
simple example requires a strengthening of our definition of the a.s.~
convergence of $\widehat{\Theta}_{n}$. This particular problem seemed to
have been implicitly identified by \citet{Ziezold1977}, as this author
proposed the following type of convergence, which specializes the event
presented in equation (\ref{eq:as definition2}),
\begin{equation}
    \lt\lb \omega\in\Omega: \bigcap_{n=1}^{\infty}\overline{\bigcup_{m=n}^{\infty}
    \widehat{\Theta}_{m}(\omega)} \subseteq \Theta \rt\rb,
   \label{eq:as definition3}
\end{equation}
where $\overline{A}$ indicates the \textit{closure} of set $A$ in $\cX$.
For convenience, this particular type of convergence will be denoted by
$\limcsup\,\widehat{\Theta}_{n}\subseteq\Theta$, a.s., where the $\climsup$ 
operator is here defined with respect to set inclusion on the power set of $\cX$.
It is easy to see why definition (\ref{eq:as definition3}) resolves the
issue illustrated in example \ref{exa:closure}. By taking the closure
of $\bigcup_{m=n}^{\infty}\widehat{\Theta}_{m}$, we include all the
elements for which there exists a sequence of $\hat{\theta}_{n}$'s
converging to $\E[X]$, and therefore for real-valued random variables, 
\begin{equation} 
    \E[X] \in \overline{\bigcup_{m=n}^{\infty}\widehat{\Theta}_{m}},
   \notag
\end{equation}
for every $n$, which implies that $\limcsup\,\widehat{\Theta}_{n}=\{\E[X]\}$, as
desired, thereby ensuring complete agreement between the classical
and Fr\'{e}chet inferential approaches for this particular example.
Note that these issues are neither related to the completeness of the
underlying space of interest, nor associated to the question of the non-emptiness of
$\Theta$. 

Since \citet{Sverdrup1981} assumed that $\cX$ is compact, it follows
that $\Theta$ and $\widehat{\Theta}_{n}$ are non-empty, in this
case. The separability of $\cX$ is not sufficient to ensure that 
$\Theta$ and the $\widehat{\Theta}_{n}$'s are non-empty. 
Nonetheless, observe that if $\,\widehat{\Theta}_{n}=\varnothing$, then the events
in equations (\ref{eq:as definition2}) and (\ref{eq:as definition3}) are
trivially almost certain, since $\varnothing\subseteq A$, for all
$A\subseteq\cX$.

\subsection{Kuratowski Upper Limit}\label{sec:kuratowski}
It can easily be shown that the type of convergence envisaged
by \citet{Ziezold1977} is, in fact, equivalent to the celebrated 
upper limit introduced by \citet{Kuratowski1966}, which has been
adopted as the preferred type of convergence in set-valued analysis 
\citep[see][]{Aubin1990}. The Kuratowski upper
limit is defined over a metric space $(\cX,d)$, for some sequence of subsets
$A_{n}\subseteq\cX$, as follows
\begin{equation}
  \begin{aligned}
    \Limsup_{n\to\infty} A_{n}
    &:= \lt\lb x\in\cX : \climinf_{n\to\infty} d(x,A_{n})=0 \rt\rb \\
    &= \Big\lb x\in\cX : \big\{A_{n}\cap
    N_{\epsilon}(x)\neq\varnothing\big\} \te{ i.o.}, \for\epsilon>0\Big\rb,
    \label{eq:Limsup}
  \end{aligned}
\end{equation}
where $\climinf$ and $\Limsup$ are taken with respect to real numbers and
subsets of $\cX$, respectively, and with
$N_{\epsilon}(x):=\{x\pri\in\cX:d(x,x\pri)<\epsilon\}$.
The second formulation of $\Limsup$ in equation (\ref{eq:Limsup})
immediately follows from the positivity of the metric. Also, observe
that the Kuratowski upper limit is equivalent to the set of
\textit{cluster points} of the sequences, $x_{n}\in A_{n}$ \citep{Aubin1990}.
Clearly, the Kuratowski upper limit of any sequence of sets is
\textit{closed}, and moreover, it contains the conventional
set-theoretical upper limit, such that for any sequence of random
sets $A_{n}$, 
\begin{equation}
      \climsup_{n\to\infty} A_{n} \subseteq \Limsup_{n\to\infty} A_{n}.
      \notag
\end{equation}
Importantly, it can be easily shown that the Kuratowski upper limit and the
quantity studied by \citet{Ziezold1977} are equivalent, as
stated in the following lemma. 
\begin{lem}\label{lem:kuratowski}   
  Given a metric space $(\cX,d)$, for any sequence of sets
  $A_{n}\subseteq\cX$, 
  \begin{equation}
       \limcsup_{n\to\infty} A_{n} = \Limsup_{n\to\infty} A_{n}.
       \notag
  \end{equation}
\end{lem}
\begin{proof}
  Clearly, $\limcsup A_{n} = \varnothing$, if and only
  if, $\Limsup A_{n} = \varnothing$. Thus, assume that these two outer
  limits are non-empty, and choose $x_{0}\in\limcsup A_{n}$. Then,
  $x_{0}\in\overline{\bigcup_{m=N}^{\infty}A_{m}}$ for every $N$ and 
  there exists a subsequence $x_{k}$ such that $x_{k}\in A_{n_{k}}$,
  for every $k$, which satisfies $x_{k}\to x_{0}$.  
  Hence, we have $\climinf d(x_{0},A_{n})=0$, and by definition
  (\ref{eq:Limsup}), $\limcsup A_{n} \subseteq \Limsup A_{n}$. 

  Conversely, choose $x_{0}\in\Limsup A_{n}$. Then, there exists a
  subsequence $x_{k}$ such that $x_{k}\in A_{n_{k}}\cap
  N_{\epsilon}(x_{0})$, for every $k$ and for every $\epsilon>0$, which
  satisfies $x_{k}\to x_{0}$, as $k\to\infty$. 
  This implies that $x_{0}\in
  \bigcap_{N=1}^{\infty}\overline{\bigcup_{m=N}^{\infty}A_{m}}$, 
  and therefore $\limcsup A_{n} \supseteq \Limsup A_{n}$, which
  completes the proof. 
\end{proof}
Observe that $\Limsup A_{n}$ can be empty. Consider the following
diverging sequence of sets, $A_{n}:=[n-1,n+1]$, for every $n$. It is
immediate that $\Limsup A_{n}=\varnothing$. 
Throughout the rest of the paper, we will neither assume the existence
nor the uniqueness of $\Theta^{r}$ and the $\Theta^{r}_{n}$'s. In
particular, in the sequel, $\Theta^{r}$ may be empty, a subset of
$\cX$, or a singleton set.

\section{Almost Sure Consistency of Fr\'{e}chet Sample
  Mean}\label{sec:main result}
In this section, we prove a strong law of large numbers 
for sample Fr\'echet means in spaces having a bounded
metric. This result can be regarded as an adaptation of Ziezold's
(1977) original result to spaces equipped with a bounded metric. 
This new proof also allows us to re-formulate Ziezold's theorem
using the Kuratowski upper limit. 
\begin{thm}\label{thm:main}   
   Given a probability space $(\Omega,\cF,\p)$ and a separable
   bounded metric space $(\cX,d)$, let $X_{1},\ldots,X_{n}$ be a sequence of
   independent and identically distributed (iid)
   abstract-valued random variables, such that
   $X_{i}:\Omega\mapsto\cX$, for every $X_{i}$.
   Then, 
   \begin{equation}\label{eq:main}
         \widehat{\sig}^{r}_{n} \to \sig^{r}
         \quad\te{a.s.,}
         \quad\quad
         \te{and}
         \quad\quad
         \Limsup_{n\to\infty}\,\widehat{\Theta}^{r}_{n}\subseteq\Theta^{r}
         \quad\te{a.s.,}
         \notag
   \end{equation}
   for every finite $r\geq1$, 
   and where $\Limsup$ is defined as in equation (\ref{eq:Limsup}). 
\end{thm}
The particular mode of convergence of the Fr\'{e}chet sample mean used in
theorem \ref{thm:main} will sometimes be denoted by
$X_{n}\stack{a.s.}{\to} X$, which implies that 
$\Limsup X_{n} \subseteq X$ with probability one. 
Observe that the integrability of the $r\tth$ order metric is implied by the
finiteness of both $d$ and $\mu$. Since $d(x,y)\leq M$, for every $x,y\in\cX$,
we have for any arbitrary $\alpha\in\cX$ and finite $r\geq1$,
\begin{equation}
       \E[d(X,\alpha)^{r}] = \int_{\cX} |d(x,\alpha)|^{r}d\mu(x) \leq 
       \int_{\cX}M^{r}d\mu(x) = M^{r}\mu(\cX) < \infty,
        \notag
\end{equation}
by the linearity of the Lebesgue integral, and the fact that $\mu$
is a probability measure. 
The integrability of the exponentiated metric was not
explicitly assumed by \citet{Sverdrup1981}. This author, however,
assumed that $\cX$ is compact, which implies that $d^{r}$ is
integrable for any finite $r\geq1$. 

The key to the proof of theorem \ref{thm:main} is based on a classical result, due to
\citet{Rao1962}, which stipulates the conditions under which the weak
convergence of a probability measure is equivalent
to the uniform convergence of a probability measure, in a sense made clear
in theorem \ref{thm:rao}. This can be seen as a generalization of the 
Glivenko-Cantelli lemma to random variables taking values in
separable metric spaces \citep[see also][chap.~ 2]{Parthasarathy1967}. In this
result, we will need to define a class of functions
on the separable space $\cX$, which we will denote by $\cF:=\cF(\cX)$,
whereby every $f\in\cF$ is a real-valued continuous function that
satisfies $f:\cX\mapsto\R$. Such a class of functions is said to be
\textit{uniformly bounded} when for every $f\in\cF$, and every $x\in\cX$,
there exists an $M\in\R$, such that $f(x)\leq M$. In addition,
$\cF$ is \textit{equicontinuous at a point} $x_{0}\in\cX$, if for
every $\epsilon>0$, there exists $\delta(x_{0})>0$, such that for every $u\in
N_{\delta}(x_{0}):=\{u\in\cX:d(x_{0},u)<\delta\}$, we have
$|f(x)-f(u)|<\epsilon$, for every $f\in\cF$. The class $\cF$ is said
to be equicontinuous if it is equicontinuous for every
$x\in\cX$. Finally, $\cF$ is said to be \textit{uniformly
  equicontinuous} if $\delta$ does not depend on $x_{0}$.
We will denote the collection of all finite measures on $\cB$ by
$\cM(\cB)$, and $\Rightarrow$ will indicate weak convergence. 
\begin{thm}[Rao, 1962, p.672]\label{thm:rao}
  Let $\cF(\cX)$ be a class of real-valued functions on a separable
  space $\cX$, and assume that $\cF(\cX)$ is (i) dominated by a continuous
  integrable function on $\cX$, and that (ii) $\cF(\cX)$ is
  equicontinuous. If, for some sequence of measures
  $\mu_{n}\in\cM(\cB)$, and $\mu\in\cM(\cB)$, we have
  $\mu_{n}\Rightarrow\mu$, a.s., then
  \begin{equation}
      \lim_{n\to\infty}\sup_{f\in\cF} \lt|\int fd\mu_{n} -
      \int fd\mu\rt|=0, \quad\te{a.s..}
        \notag
  \end{equation}
\end{thm}

The following lemma will be used in the proof of theorem
\ref{thm:main}. This result links the properties of a bounded metric
space with the conditions required in Rao's (1962) theorem. For this
purpose, we will require the following classes of \textit{point functions} on a metric
space \citep[see][]{Searcoid2007}.
\begin{dfn}\label{dfn:point}
   For any metric space $(\cX,d)$, the $z$-point function is defined
   as $d_{z}(x):=d(z,x)$ for every $x\in\cX$. The class of point functions
   on $(\cX,d)$ is then denoted by $\cD(\cX) := \lb d_{z}: \for
   z\in\cX\rb$. Similarly, we will make use of the class of exponentiated
   point functions, defined as follows, 
   \begin{equation}
        \cD^{r}(\cX) :=  \lt\lb d^{r}_{z}: \for z\in\cX\rt\rb,
        \notag
   \end{equation}
   for every finite $r\geq1$, and 
   where elements in either $\cD$ or $\cD^{r}$ will be denoted by
   $d_{z}$, and $d_{z}^{r}$, respectively.
\end{dfn}
\begin{lem}\label{lem:finite}
  If $(\cX,d)$ is a bounded metric space, then $\cD^{r}(\cX)$ is
  uniformly bounded and uniformly equicontinuous for every finite $r\geq1$.
\end{lem}
\begin{proof}
  By the boundedness of $(\cX,d)$, there exists an $M\in\R$, such that 
  $d(x,y)\leq M$, for every $x,y\in\cX$. Therefore, $d_{z}(x)\leq M$,
  for every $x\in\cX$, for every $d_{z}\in\cD$, and thus
  $\cD$ is uniformly bounded. Moreover, since $d^{r}_{z}(x)\leq
  M^{r}<\infty$, for every finite $r\geq1$, it follows that each $\cD^{r}$ also
  forms a uniformly bounded class of functions. Next, by the reverse
  triangle inequality, we
  have $\lt|d_{z}(x) - d_{z}(x_{0}) \rt| \leq d(x,x_{0})$, for all
  $x,x_{0},z\in\cX$, thereby proving the (uniform) equicontinuity of
  the class $\cD$ on
  $\cX$. For the case of $r\geq1$, we consider the exponentiated
  version of the triangle inequality. Using the binomial expansion,
  \begin{equation}
    \begin{aligned}
      d(z,x)^{r} &\leq 
       \Big(d(z,x_{0}) + d(x_{0},x)\Big)^{r} \\
      & =d(z,x_{0})^{r} +
      \sum_{k=1}^{r-1}\binom{r}{k}d(z,x_{0})^{r-k}d(x_{0},x)^{k} +
      d(x_{0},x)^{r}.
    \end{aligned}
      \notag
  \end{equation}
  Similarly, for any given $x_{0}\in\cX$, 
      $d(z,x_{0})^{r} \leq d(z,x)^{r} +
      \sum_{k=1}^{r-1}\binom{r}{k}d(z,x)^{r-k}d(x,x_{0})^{k} + d(x,x_{0})^{r}$. 
  Combining these two inequalities and invoking the symmetry of $d$, we have
  \begin{equation}
    \begin{aligned}
      \lt|d(z,x)^{r} - d(z,x_{0})^{r} \rt| & \leq d(x_{0},x)^{r} +
               d(x_{0},x)M^{r-1}\sum_{k=1}^{r-1}\binom{r}{k}\\  
           & \leq d(x_{0},x)M^{r-1}\lt(1+\sum_{k=1}^{r-1}\binom{r}{k}\rt),
    \end{aligned}
    \notag
  \end{equation}
  where $M$ is the uniform bound on the class $\cD$. Now, choose
  $\delta=\epsilon/\gamma M^{r-1}$, where $\gamma:=1+\sum_{k=1}^{r-1}\binom{r}{k}$,
  such that if $d(x,x_{0})<\delta$, then
  $\lt|d^{r}_{z}(x) - d^{r}_{z}(x_{0}) \rt| < \gamma\delta M^{r-1} = \epsilon$,
  for every $x\in N_{\delta}(x_{0})$, for every $d_{z}^{r}\in\cD^{r}$, 
  thence proving the equicontinuity of $\cD^{r}$ at $x_{0}$. Since $\delta$ did
  not depend on the choice of $x_{0}$, it follows that $\cD^{r}$ is
  also uniformly equicontinuous. 
\end{proof}
\begin{proof}[Proof of Theorem \ref{thm:main}]
    Observe that the theorem is trivially
    verified if $\Limsup\widehat{\Theta}^{r}_{n}=\varnothing$. Thus, 
    assume that $\Limsup\widehat{\Theta}^{r}_{n}$ is non-empty.
    We here adopt the line of argument followed by 
    \citet{Sverdrup1981}. However, since we are not assuming
    compactness, there are several aspects of Sverdrup-Thygeson's proof that
    becomes somewhat delicate. In the sequel, we will make use of the following
    quantities formulated with respect to the class of point functions
    described in definition \ref{dfn:point}. For every $z\in\cX$, let
    \begin{equation}
          T_{n}(z) := 
          \frac{1}{n}\sum_{i=1}^{n} d^{r}_{z}(X_{i}) -
          \int_{\cX} d^{r}_{z}(x)d\mu(x), 
          \label{eq:Tn}
    \end{equation}
    and similarly, 
    \begin{equation}
          T\as_{n}(z) := 
          \frac{1}{n}\sum_{i=1}^{n} d^{r}_{z}(X_{i}) -
          \int_{\cX} d^{r}_{\theta}(x)d\mu(x).
          \label{eq:Tnas}
    \end{equation}
    Since $T_{n}(x)$ is real-valued, one can invoke the 
    strong law of large numbers for real-valued random
    variables, which gives
    \begin{equation}
        T_{n}(z) \to 0, \quad \te{a.s.}, \quad\for z\in\cX.
        \label{eq:as}
    \end{equation}

    Note, however, that since we have used
    infima in the definitions of the Fr\'{e}chet theoretical and sample
    means in equations (\ref{eq:frechet theory}) and 
    (\ref{eq:frechet sample}), it follows that the convergence of
    $T_{n}(z) \to 0$ is not assured when $z$ is an element of $\Theta$
    or an element of $\widehat{\Theta}_{n}$. However, as established in
    lemma \ref{lem:finite}, the class of point functions, $\cD^{r}(\cX)$, is
    uniformly bounded and (uniformly) equicontinuous. Moreover, we
    have seen that the finiteness of $\E[d^{r}_{z}(X)]$ is implied by
    the boundedness of $d$, such that $\E[d^{r}_{z}(X)]\leq M^{r}\mu(\cX)$. Thus, it
    follows that there exists a continuous 
    integrable function, i.e.~ $f(x):=M^{r}$, dominating every
    $d_{z}^{r}\in\cD^{r}$. Moreover, a classical result on the
    convergence of empirical measures based on iid random variables
    taking values in separable metric spaces
    \citep[see][theorem 7.1, p.53]{Parthasarathy1967} implies that 
    \begin{equation}
        \mu_{n}\Rightarrow \mu, \quad\te{a.s.,} 
        \label{eq:weak as}
    \end{equation}
    where $\mu_{n}:=n^{-1}\sum_{i=1}^{n}\delta_{X_{i}}$, is the
    empirical measure on $\cX$. Therefore, we are in a position
    to apply theorem \ref{thm:rao}, which shows that the empirical
    measure, $\mu_{n}$, converges uniformly with probability 1. That is, 
    \begin{equation}
        \p\lt[\sup_{z\in\cD^{r}}
        \lt|\frac{1}{n}\sum_{i=1}^{n} d^{r}_{z}(X_{i}) -
          \int_{\cX}  d^{r}_{z}(x)d\mu(x)\rt|     
        \to 0\rt]=1, 
        \notag
    \end{equation}
    which may be re-written as     
    \begin{equation}
          \sup_{z\in\cD^{r}}\big|T_{n}(z)\big| =
          \sup_{z\in\cX}\big|T_{n}(z)\big| \to 0, \quad\te{a.s.}.
          \label{eq:unif as}
    \end{equation}
    Consequently, $T_{n}(\hat{\theta}_{n})\to 0$, a.s., and
    $T_{n}(\theta)\to 0$, a.s., for every
    $\hat{\theta}_{n}\in\widehat{\Theta}_{n}$ and every $\theta\in\Theta$,
    respectively. 

    Further, from the definition of $\hat{\theta}_{n}$ and $\theta$, we can
    `sandwich' $T_{n}\as(\hat{\theta}_{n})$ in the following
    manner. Firstly, observe that by the minimality of the 
    $\theta$'s, 
    \begin{equation}
      \begin{aligned}
        T_{n}(\hat{\theta}_{n}) &= \frac{1}{n}\sum_{i=1}^{n}
        d^{r}_{\hat{\theta}_{n}}\!\!(X_{i}) - \int_{\cX}
        d^{r}_{\hat{\theta}_{n}}\!\!(x)d\mu(x) \\ 
        &\leq \frac{1}{n}\sum_{i=1}^{n} d^{r}_{\hat{\theta}_{n}}\!\!(X_{i}) -
        \int_{\cX} d^{r}_{\theta}(x)d\mu(x) = T_{n}\as(\hat{\theta}_{n}).
      \end{aligned}
      \label{eq:sandwich1}
    \end{equation}
    Secondly, by the minimality of the $\hat{\theta}_{n}$'s, we similarly have,
    \begin{equation}
      \begin{aligned}
        T_{n}\as(\hat{\theta}_{n}) &= \frac{1}{n}\sum_{i=1}^{n}
        d^{r}_{\hat{\theta}_{n}}\!\!(X_{i}) -
        \int_{\cX}  d^{r}_{\theta}(x)d\mu(x)\\
        &\leq \frac{1}{n}\sum_{i=1}^{n} d^{r}_{\theta}(X_{i}) - \int_{\cX}
        d^{r}_{\theta}(x)d\mu(x) = T_{n}(\theta).
      \end{aligned}
      \label{eq:sandwich2}
    \end{equation}
    Thence, combining equations (\ref{eq:sandwich1}) and (\ref{eq:sandwich2}),
    we obtain, 
    \begin{equation}
              T_{n}(\hat{\theta}_{n}) 
              \leq 
              T_{n}\as(\hat{\theta}_{n})
              \leq 
              T_{n}(\theta),
        \notag
    \end{equation}
    such that, using equation (\ref{eq:unif as}),
    \begin{equation}
        |T_{n}\as(\hat{\theta}_{n})| \leq
        \max\{|T_{n}(\hat{\theta}_{n})|,|T_{n}(\theta)|\} \to 0,
        \quad\te{a.s.,}
        \label{eq:sigma convergence}
    \end{equation}
    which proves the a.s.~ convergence of $\widehat\sig^{r}_{n}$ to $\sig^{r}$.

    We now turn to the convergence properties of the Fr\'{e}chet sample
    mean of the $r\tth$ order, $\widehat{\Theta}^{r}_{n}$. Here, we
    generalize Ziezold's (1977) proof strategy to Fr\'{e}chet sample means of
    any order \citep[see also][p.185]{Molchanov2005}. Choosing
    \begin{equation}
         \hat{\theta}\in\Limsup_{n\to\infty}\,\widehat{\Theta}^{r}_{n},
         \notag
    \end{equation}
    it then suffices to show that $\hat{\theta}\in\Theta^{r}$, which
    is verified if $\E[d(X,\hat\theta)^{r}]\leq\E[d(X,x\pri)^{r}]$,
    for every $x\pri\in\cX$. We proceed by constructing the following
    subsequence of natural numbers. 

    Observe that from the definition of the Kuratowski upper limit
    and the equivalence relation reported in lemma
    \ref{lem:kuratowski}, it follows that
    $\hat{\theta}\in\op{Cl}(\bigcup_{m=n}^{\infty}\widehat{\Theta}_{m}^{r})$,
    for every $n$, where $\op{Cl}(\cdot)$ denotes the closure of a
    set. Thus, one can construct a subsequence, $\{n_{k}:k\in\N\}$, such
    that for every $k$, there exists an element
    $\hat{\theta}_{k}\in\bigcup_{m=k}^{\infty}\widehat{\Theta}^{r}_{m}$,
    which satisfies $d(\hat{\theta}_{k},\hat{\theta})\leq
    1/k$. Moreover, we can define 
    $n_{k}:=\min\{n\in\Na: n\geq
    k,\hat{\theta}_{k}\in\widehat{\Theta}^{r}_{n}\}$. Now, after
    an application of the triangle inequality, followed by the
    Minkowski inequality, we have
    \begin{equation}
        \lt(\frac{1}{n_{k}}\sum_{i=1}^{n_{k}}d(X_{i},\hat{\theta})^{r}\rt)^{1/r}
        \leq
        \lt(\frac{1}{n_{k}}\sum_{i=1}^{n_{k}}d(X_{i},\hat{\theta}_{k})^{r}\rt)^{1/r} +
        \lt(\frac{1}{n_{k}}\sum_{i=1}^{n_{k}}d(\hat{\theta}_{k},\hat{\theta})^{r}\rt)^{1/r},
        \notag
    \end{equation}
    which gives 
    \begin{equation}
        \lt(\frac{1}{n_{k}}\sum_{i=1}^{n_{k}}d(X_{i},\hat{\theta})^{r}\rt)^{1/r}
        \leq \lt(\frac{1}{n_{k}}\sum_{i=1}^{n_{k}}d(X_{i},\hat{\theta}_{k})^{r}\rt)^{1/r} +
        \frac{1}{k}. 
        \notag
    \end{equation}
    As $k\to\infty$, it then follows from equation (\ref{eq:unif as})
    that since $(n_{k})_{k\in\Na}$ is a subsequence of
    $(n)_{n\in\Na}$, we obtain
    \begin{equation}
        \lt(\E[d(X,\hat{\theta})^{r}]\rt)^{1/r} 
        \leq
        \climinf_{k\to\infty}\lt(\frac{1}{n_{k}}\sum_{i=1}^{n_{k}}d(X_{i},\hat{\theta}_{k})^{r}\rt)^{1/r},
        \label{eq:ziezold1}
    \end{equation}
    where $\climinf$ is here taken with respect to non-negative real numbers.
    Moreover, by construction, each $\hat{\theta}_{k}$ is minimal with
    respect to any element $x\pri\in\cX$, such that 
    \begin{equation}        
        \frac{1}{n_{k}}\sum_{i=1}^{n_{k}}d(X_{i},\hat{\theta}_{k})^{r}
        \leq
        \frac{1}{n_{k}}\sum_{i=1}^{n_{k}}d(X_{i},x\pri)^{r},
        \label{eq:ziezold2}
    \end{equation}
    for every $x\pri\in\cX$ and $k\in\Na$. Observe that given the
    continuity and monotonicity of $g(x):=x^{1/r}$ on positive real numbers, we
    have $\climinf g(x_{n})=g(\climinf x_{n})$, for every sequence
    satisfying $x_{n}\in\R^{+}$. Therefore, it suffices to
    combine equations (\ref{eq:ziezold1}) and (\ref{eq:ziezold2}) in
    order to obtain $\E[d(X,\hat{\theta})^{r}] \leq \E[d(X,x\pri)^{r}]$,
    for every $x\pri\in\cX$, as required. Thence, $\hat{\theta}\in\Theta^{r}$ a.s.,
    but since $\hat{\theta}$ was arbitrary, we have
    $\Limsup\,\widehat{\Theta}^{r}_{n}\subseteq\Theta^{r}$ a.s., as
    required.
\end{proof}

\section{Restricted Fr\'{e}chet Means}\label{sec:restricted}
Theorem \ref{thm:main} can be extended to the case of the restricted
Fr\'{e}chet mean. This is a concept that was originally introduced and
studied by \citet{Sverdrup1981}. Interest in restricted Fr\'{e}chet means is
motivated by the fact that the domain of some abstract-valued
random variables may be too large to be optimized in a reasonable 
amount of time. This perspective is especially relevant when
considering discrete metric spaces of graphs, where minimization may be
computationally NP-hard.

In such cases, the Fr\'{e}chet sample mean may be more
suitably defined as one of the elements in the sample at hand. That is,
consider the following definition of the \textit{restricted} Fr\'{e}chet sample
mean and variance,
\begin{equation}
      \widehat{\Theta}^{\ast,r}_{n}:=\argmin_{x\pri\in \bX}
      \sum_{i=1}^{n}d(X_{i},x\pri)^{r} 
      \quad\te{and}\quad
      \widehat{\sig}_{n}^{\ast,r} :=
      \min_{x\pri\in\bX}\sum_{i=1}^{n}d(X_{i},x\pri)^{r},
      \notag
\end{equation}
where $\bX:=\{X_{1},\ldots,X_{n}\}\subseteq \cX$ denotes the set of sampled
variables. In practice, the sample mean is chosen
among the available sampled iid realizations from $X$. In particular,
observe that we have employed the minimum instead of the infimum in the
definitions of both $\widehat{\Theta}^{\ast,r}_{n}$ and
$\widehat{\sig}^{\ast,r}_{n}$, as the required optimal values 
necessarily exist, albeit they may not be unique. Hence,
observe that $\widehat{\Theta}^{\ast,r}_{n}\neq\varnothing$ for any $n$.
Theoretical analogues of these restricted quantities can be defined as
follows, 
\begin{equation}
     \Theta^{\ast,r} := \argmin_{x\pri\in W} \int\limits_{\cX}d(x,x\pri)^{r}
     d\mu(x), 
     \quad\te{and}\quad
     \sig^{\ast,r} := \min_{x\pri\in W}\int\limits_{\cX}d(x,x\pri)^{r}d\mu(x),
      \notag
\end{equation}
where $W$ is the support of $\mu$, denoted $\op{supp}(\mu)$, and is
assumed to be \textit{closed}. Observe that this closure condition is
required in order to ensure that the Fr\'{e}chet mean is contained within $\op{supp}(\mu)$.
As previously, the elements of $\Theta\as$ and
$\widehat{\Theta}\as_{n}$ will be denoted by $\theta\as$'s and
$\hat{\theta}\as_{n}$'s, respectively. We here prove a generalization of a consistency
result due to \citet{Sverdrup1981} on the a.s.~ convergences of the
restricted Fr\'{e}chet sample mean and variance. 
\begin{thm}\label{thm:restricted}
   Under the conditions of theorem \ref{thm:main}, for every $r\geq1$,
   and assuming that $\op{supp}(\mu)$ is closed, 
   \begin{equation}\label{eq:main}
         \widehat{\sig}^{\ast,r}_{n} \to \sig^{\ast,r}
         \quad\te{a.s.,}
         \quad\quad
         \te{and}
         \quad\quad
         \Limsup_{n\to\infty}\,\widehat{\Theta}^{\ast,r}_{n} \subseteq \Theta^{\ast,r}
         \quad\te{a.s..}
      \notag
   \end{equation}  
\end{thm}
\begin{proof} 
   Let us denote a quantity analogous to the ones defined in equations
   (\ref{eq:Tn}) and (\ref{eq:Tnas}), but here based on the
   \textit{restricted} theoretical Fr\'echet mean,
   \begin{equation}
          \op{TR}\as_{n}(z) := 
          \frac{1}{n}\sum_{i=1}^{n} d^{r}_{z}(X_{i}) -
          \int_{\cX} d^{r}_{\theta\as}(x)d\mu(x),
          \label{eq:TRnas}
   \end{equation}
   where $\theta\as\in\Theta\as$. We will first demonstrate that 
   \begin{equation}
        \min_{x\pri\in\bX} \Big|\op{TR}_{n}\as(x\pri) -
        \op{TR}\as_{n}(\theta\as)\Big|\to 0, \quad\te{a.s..}
        \label{eq:TR as}
   \end{equation}
   In order to prove this a.s.~ convergence, we need the following
   quantity, 
   \begin{equation}
         s(\delta):= \sup_{z\in W}\,\sup_{d(x,y)<\delta}\big| 
         d^{r}_{z}(x) - d^{r}_{z}(y)\big|, 
        \label{eq:s(delta)}
   \end{equation}
   where the second supremum is taken over all pairs of elements
   $x,y\in W$, satisfying $d(x,y)<\delta$. Since the class of
   exponentiated point functions on $\cX$, denoted
   $\cD^{r}$, was shown to be uniformly equicontinuous in lemma
   \ref{lem:finite}, it follows that $s(\delta)\to0$, as
   $\delta\to0$. Moreover, it is straightforward
   to see that for every $\delta>0$, we have
   \begin{equation}
     \begin{aligned}
       \sup_{d(x,y)<\delta}\big|\op{TR}_{n}\as(x)-\op{TR}_{n}\as(y)\big|
       &= \sup_{d(x,y)<\delta}\lt|
       \frac{1}{n}\sum_{i=1}^{n}d^{r}_{x}(X_{i}) -
       \frac{1}{n}\sum_{i=1}^{n}d^{r}_{y}(X_{i})\rt|\\
       &\leq \sup_{d(x,y)<\delta}
       \frac{1}{n}\sum_{i=1}^{n}\Big|d^{r}_{x}(X_{i}) -
       d^{r}_{y}(X_{i})\Big|\\
       &\leq s(\delta).
      \notag
     \end{aligned}
   \end{equation}
   Next, let $O_{\delta}:=\{x\in\cX:d(x,\theta\as)<\delta\}$, for any
   $\delta>0$. Since $\theta\as\in\op{supp}(\mu)$, from the definition
   of the restricted Fr\'{e}chet mean, it follows that
   $\mu(O_{\delta})=:\alpha>0$. Hence, 
   \begin{equation}
         \p\lt[\lb X_{1}\in O_{\delta}\rb\cup\hdots\cup \lb X_{n}\in O_{\delta}\rb\rt]
         =1-\prod^{n}_{i=1}\p\lt[\lb X_{i}\notin O_{\delta}\rb\rt] = 1-(1-\alpha)^{n},
      \notag
   \end{equation}
   which converges to $1$, as $n\to\infty$, for any
   $\alpha>0$. Moreover, observe that since $x\pri\in W$, for every
   $x\pri\in\bX$, we also have 
   \begin{equation}
     \begin{aligned}
       \climsup_{n\to\infty}\min_{x\pri\in\bX}\big|\op{TR}_{n}\as(x\pri)-\op{TR}_{n}\as(\theta\as)\big|
       & \leq s(\delta).
      \notag
     \end{aligned}
   \end{equation}
   It then suffices to let $\delta\to0$, in order to obtain equation
   (\ref{eq:TR as}). Now, from the definitions of $\op{TR}_{n}\as$ and
   $T_{n}$, it can be seen that
   $\op{TR}_{n}\as(\theta\as)=T_{n}(\theta\as)$, and therefore 
   \begin{equation}
       \op{TR}\as_{n}(\hat{\theta}\as_{n}) = \min_{x\pri\in\bX} \op{TR}\as_{n}(x\pri) \leq T_{n}(\theta\as)
        + \min_{x\pri\in\bX} \big|\!\op{TR}\as_{n}(x\pri) - \op{TR}\as_{n}(\theta\as)\big|,
        \notag
   \end{equation}
   by the optimality of $\hat{\theta}\as_{n}$.
   This can be bounded below by using the minimality of $\theta\as$, such that 
   \begin{equation}
     \begin{aligned}
       T_{n}(\hat{\theta}\as_{n}) &= \frac{1}{n}\sum_{i=1}^{n}
       d^{r}_{\hat{\theta}\as_{n}}\!\!(X_{i}) - \int_{\cX}
       d^{r}_{\hat{\theta}\as_{n}}\!\!(x)d\mu(x) \\
       & \leq \frac{1}{n}\sum_{i=1}^{n} d^{r}_{\hat{\theta}\as_{n}}\!(X_{i}) -
       \int_{\cX} d^{r}_{\theta\as}(x)d\mu(x) =\op{TR}\as_{n}(\hat{\theta}\as_{n}).
     \end{aligned}
     \notag
   \end{equation}
   Combining the last two results, we obtain the following `sandwich'
   inequality of $\op{TR}\as_{n}(\hat{\theta}\as_{n})$, 
   \begin{equation}
       T_{n}(\hat{\theta}\as_{n}) \leq
       \op{TR}\as_{n}(\hat{\theta}\as_{n}) \leq  T_{n}(\theta\as)
        + \min_{x\pri\in\bX} \big|\!\op{TR}\as_{n}(x\pri) - \op{TR}\as_{n}(\theta\as)\big|.
       \notag
   \end{equation}
   Thence, this gives a.s.,
   \begin{equation}
      |\op{TR}\as_{n}(\hat{\theta}\as_{n})| \leq 
      \max\big\{|T_{n}(\hat{\theta}\as_{n})|, 
      |T_{n}(\theta\as)|
        + \min_{x\pri\in\bX} \big|\!\op{TR}\as_{n}(x\pri) - \op{TR}\as_{n}(\theta\as)\big|
      \big\} \to 0, 
      \notag
   \end{equation}
   using the strong law of large numbers on $T_{n}(\theta\as)$, and using
   equation (\ref{eq:TR as}) for the second term in the maximum. This
   proves that $\widehat{\sig}_{n}\to\sig$, a.s.. 
   The proof of $\Limsup\,\widehat{\Theta}\as_{n} \subseteq
   \Theta\as$ with probability 1, can be conducted using the same
   construction described in the proof of theorem \ref{thm:main}, by choosing
   $\hat{\theta}\as\in\Limsup\,\widehat{\Theta}_{n}\as$, and noting
   that $\op{supp}(\mu)$ was assumed to be closed.
\end{proof}
\begin{rmk}\label{rmk:equicontinuity}
  The use of uniform equicontinuity in the proof of theorem
  \ref{thm:restricted} requires special mention. \citet{Sverdrup1981}
  was able to invoke the continuity of $s(\delta)$ with respect to
  $\delta$ in equation (\ref{eq:s(delta)}) by using the compactness of
  $\cX$. Here, this property immediately follows from the uniform
  equicontinuity of the class of exponentiated point functions,
  $\cD^{r}(\cX)$. This was the sole argument in the proof of
  \citet{Sverdrup1981} for the a.s.~ convergence of the restricted Fr\'{e}chet
  sample mean that required the compactness of
  $\cX$. Hence, the boundedness of $d$ constitutes a sufficient
  condition.
\end{rmk}
\begin{rmk}\label{rmk:support}
   Under our assumptions and the ones postulated by both 
   \citet{Ziezold1977} and \citet{Sverdrup1981}, there is no guarantee
   that $\Theta\subseteq\op{supp}(X)$ holds, as assumed in the definition of
   the restricted Fr\'{e}chet mean. In particular, one can easily
   construct a measure space where $\Theta$ belongs to a set of $\mu$-measure
   zero. Consider the random variable described in example
   \ref{exa:closure}, where two point masses were located at $-1$ and
   $1$, respectively, and the Fr\'{e}chet mean was computed with respect
   to the square of the Manhattan distance. Clearly,
   the Fr\'{e}chet mean is located in the barycenter of the interval $[-1,1]$
   but that center of mass does not belong to $\op{supp}(X)$, which is
   simply $\{-1,1\}$.
\end{rmk}

\section{Conclusion}\label{sec:conclusion}
In this paper, we have generalized the results due to
\citet{Sverdrup1981} by relaxing the compactness assumption made by
this author. This task has highlighted interesting links between the 
Sverdrup-Thygeson's proof and
another classical proof of the a.s.~ convergence of the Fr\'{e}chet sample
mean, due to \citet{Ziezold1977}. In particular, we have shown that
by assuming the boundedness of the metric of interest, we can deduce
the uniform boundedness and uniform equicontinuity of any family of point
functions on $\cX$. These two properties were found to be required 
on two distinct occasions when proving asymptotic convergence results
for the unrestricted
and restricted Fr\'{e}chet sample means, respectively. In the original
proof of \citet{Sverdrup1981}, these two arguments rely on
compactness, thereby showing that uniform boundedness and 
uniform equicontinuity constitute appropriate weaker assumptions.

Throughout, we have assumed that the underlying metric of interest is
a full metric. However, as was originally done by \citet{Ziezold1977},
it can be shown that our results also hold for bounded pseudo-metrics,
where one relaxes the \textit{axiom of coincidence}. 
In this case, $d(x,y)=0$ does not necessarily imply that $x=y$. 
It is easy to check that this particular property was not used in this
paper, and therefore that the aforementioned convergence theorems
remain valid for Fr\'{e}chet sample mean sets defined over separable bounded
\textit{pseudo}-metric spaces. These results may be of special interest
to statisticians considering graph-valued random variables, which are
commonly defined over bounded metric spaces. 
Future work may include the consideration of the convergence in law of
such graph-valued random variables, or concentrate on studying the
asymptotic properties of statistics summarizing the distances between
two or several groups of graphs. 

\small
\singlespacing
\bibliography{/home/cgineste/ref/bibtex/Statistics,%
             /home/cgineste/ref/bibtex/Neuroscience}
\bibliographystyle{/home/cgineste/ref/style/oupced3}

\end{document}